# Calculation of a $K_2$ group of an $F_2$ coefficients noncommutative group algebra

LiangYi Xiong[1]    GuoPing Tang[1]

(School of Mathematical Sciences, University of Chinese Academy of Sciences, Address: No.19 Yuquan Road, Shijingshan District, Beijing (A), Zip code:100049, E-mail:1435429176@qq.com; tanggp@ucas.ac.cn)

Abstract

In this paper, the $K_2$ group of $F_2$ coefficient group algebra of a noncommutative group with 8 elements(dihedral group $D_4$) is calculated, which is divided into three parts: The first part is the introduction of basic knowledge related to algebra K-theory, and a method of Magurn to calculate finite field coefficients noncommutative finite group algebra in reference [2]. In the second part, operation laws of Dennis-Stein symbols is introduced, and we combined it with the fact that $F_2[D_4]$ is a local ring to determind the direct sum term of $K_2(F_2[D_4])$ can only be $Z_2$ or $Z_4$. In the third part, we continue to make use of the fact that $F_2[D_4]$ is a local ring, and proved that the group $D_1(F_2[D_4])$ is an abelian group closely related to the group $K_2(F_2[D_4])$ through operating Dennis-Stein symbols. Then, we used group homology and the Kunneth formula of the finite abelian group version to calculate all cases of $H_2(D_1(F_2[D_4]), Z)$, and substituted the obtained results into the long exact sequence derived from the Hochschild-Serre spectral sequence for testing, and finally constructed the result: $K_2(F_2[D_4]) \cong Z_2$.

## 1. The basic knowledge of algebraic K-theory and Magurn's method of calculating noncommutative group algebra

1.1. ( Grothendieck construction) Setting $(S, \oplus)$ Is an Abelian semigroup, then the Grothendieck construction of $S$ is defined as: $G(S) \stackrel{\text{def}}{=} F(S)/< m + n - (m \oplus n) >$, in which $F(S)$ is a free Abel group generated by $S$.

Theorem 1.2. We have: $G(S) = S \times S/\sim \stackrel{\text{def}}{=} S^{-1}S$, where the equivalence relation $\sim$ is defined as: $(a,b) \sim (a',b') \Leftrightarrow \exists c \in S, s.t.: a \oplus b' \oplus c = a' \oplus b \oplus c$.

Proof: Remember that: $S \times S/\sim \overset{\text{def}}{=} S^{-1}S$, construct mapping $\varphi: G \to G(S): (m,n) \to m-n$, surjectivity is obviously, we only need to show φ is injective. Immediately, we only need to proof that:

If $m,n \in S, s.t.: [m] = [n] \in G(S)$, where [] is the canonical map of $F(S)$ to $G(S)$, then:
$\exists\ c \in S$, such that: $m \oplus c = n \oplus c: [m] - [n] = 0 \Leftrightarrow [m] - [n] = \sum([a_i + b_i] - [a_i] - [b_i]) - \sum([c_j + d_j] - [c_j] - [d_j]) \Leftrightarrow [m] + \sum([a_i] + [b_i]) + \sum[c_j + d_j] = [n] + \sum[a_i + b_i] + \sum([c_j] + [d_j])$, with lifting to the free Abel group $F(S)$, we have (Two elements in the free Abel group are equal if and only if the components of their corresponding bases are exactly the same Iin the sense of permutation order)): $m \oplus \sum(a_i \oplus b_i) \oplus \sum(c_j \oplus d_j) = n \oplus \sum(a_i \oplus b_i) \oplus \sum(c_j \oplus d_j)$, So we can Let $c = \sum(a_i \oplus b_i) \oplus \sum(c_j \oplus d_j)$. ∎

1.3. ($K_0$ Group of a ring) Let R be any ring, and Let P(R) is the set of isomorphic classes of all finitely generated projective R- modules, together with operator $\oplus$ and identity 0, forms an abelian monoid. The $K_0$ Group of ring R is defined as the group completion of P(R) (Grothendieck construction) i.e.: $K_0(R) \overset{\text{def}}{=} P^{-1}(R)P(R)$.

1.4. (General linear group) Let $R$ be any ring, the general linear group $GL_n(R)$ is the group of invertible n × n matrices. For arbitrary $g \in GL_n(R)$, by mapping $g \to \begin{pmatrix} g & 0 \\ 0 & 1 \end{pmatrix}$, $GL_n(R)$ can be embed into $GL_{n+1}(R)$, so we have: $R^\times = GL_1(R) \hookrightarrow \ldots \hookrightarrow GL_n(R) \hookrightarrow GL_{n+1}(R) \hookrightarrow \ldots$. The union of the resulting sequence $R^\times = GL_1(R) \hookrightarrow \ldots \hookrightarrow GL_n(R) \hookrightarrow GL_{n+1}(R) \hookrightarrow \ldots$ is called the (infinite) general linear group $GL(R)$.

1.5. (Elementary matrices group) Let $R$ be any ring, $i \neq j$ are distinct positive integers and $r \in R$, then the elementary matrix $e_{ij}(r) \in GL_n(R)$ which has 1 in every diagonal spot, has $r$ in the $(i,j)$-spot, and is zero elsewhere. $E_n(R) \subseteq GL_n(R)$ denotes the subgroup of $GL_n(R)$ generated by all elementary matrices $e_{ij}(r)$ with $1 \le i, j \le n$, and the union $E(R)$ of the $E_n(R)$ is the subgroup of $GL(R)$ generated by all elementary matrices. (i.e.: $E_2(R) \hookrightarrow E_3(R) \hookrightarrow \ldots \hookrightarrow E_n(R) \hookrightarrow E_{n+1}(R) \hookrightarrow \ldots$, $E(R) \overset{\text{def}}{=} \bigcup_{n=2}^{\infty} E_n(R) = \varinjlim E_n(R)$.

Remark 1.6. The commutator of two elemnetary matrices is easy to compute and simple to discribe as follows:
1. $e_{ij}(r)e_{ij}(s) = e_{ij}(r+s)$
2. $e_{ij}^{-1}(r) = e_{ij}(-r)$
3. $[e_{ij}(r), e_{kl}(s)] = \begin{cases} 1 & j \neq k, i \neq l, \\ e_{il}(rs) & j = k, i \neq l, \\ e_{kj}(-sr) & j \neq k, i = l. \end{cases}$

Definition 1.7. ($K_1$ group of a ring) Let $R$ be any ring, its $K_1$ group $K_1(R)$ is defined as: $K_1(R) \overset{\text{def}}{=} GL(R)/[GL(R), GL(R)] = GL(R)/E(R)$.

1.8. (Steinberg group of a ring) Let $R$ be any ring, the Steinberg group $St_n(R)$ of $R$ is defined as a group given by the following generators and relations:
Situation 1. $n = 2$,

Generators: $\{x_{ij}(r)|r \in R, 1 \leq i \neq j \leq 2\}$.

Generating relations: (i) $x_{ij}(r)x_{ij}(s) = x_{ij}(r+s)$,

(ii) $w_{ij}(\alpha)x_{ji}(r)w_{ij}^{-1}(\alpha) = x_{ij}(-\alpha r \alpha); w_{ij}(\alpha) \stackrel{def}{=} x_{ij}(\alpha)x_{ji}(-\alpha^{-1})x_{ij}(\alpha)$.

Situation 2. $n \geq 3$,

Generators: $\{x_{ij}(r)|r \in R, 1 \leq i \neq j \leq n\}$,

Generating relations: (ST1) $x_{ij}(r)x_{ij}(s) = x_{ij}(r+s)$

$$(\text{ST2}) \quad [x_{ij}(r), x_{kl}(s)] = \begin{cases} 1 & j \neq k, i \neq l, \\ x_{il}(rs) & j = k, i \neq l, \\ x_{kj}(-sr) & j \neq k, i = l. \end{cases}$$

There is a canonical group surjection $\phi_n: St_n(R) \to E_n(R)$, sending $x_{ij}(r)$ to $e_{ij}(r)$. So there is an obvious map $St_n(R) \to St_{n+1}(R)$. We write $St(R)$ for $St(R) \stackrel{def}{=} \varinjlim St_n(R)$.

Definition 1.9. ($K_2$ group of a ring) Let $R$ be any ring and consider the group surjection $\phi$ induced by $\phi_n$ by stabilizing: $\phi: St(R) \to E(R): x_{ij}(r) \to e_{ij}(r)$, then we define the $K_2$ group of ring $R$: $K_2(R) \stackrel{def}{=} Ker(\phi)$.

1.10. (Relative $K_2$ group of a ring) Let $R$ be any ring and $I$ be an ideal of $R$. Define the relative $K_2$ group of the ring $R$ as: $K_2(R, I) \stackrel{def}{=} Ker(K_2(R) \to K_2(R/I))$. (Note: the relative $K_2$ group of the rings here is slightly different from the definition of relative $K_2$ group of rings that in most of K-theory books. But there is a connection between them, in fact: the former is the kernel of the mapping latter to the group $K_2(R)$.

1.11. (Relative $K_1$ group of a ring) Let $R$ be any ring and $I$ be an ideal of $R$. Remember that $GL(I) = Ker(GL(R) \to GL(R/I))$, and remember $E(R, I)$ is the minimum normal subgroup of $E(R)$ includes $\{e_{ij}(x)|x \in I\}$, then we define the relative $K_1$ group as: $K_1(R, I) \stackrel{def}{=} GL(I)/E(R, I)$.

Theorem 1.12. (Long exact sequence in K-theory) Let $R$ be any ring and $I$ be an ideal of $R$. Then there is a long exact sequence: $K_2(R, I) \to K_2(R) \to K_2(R/I) \stackrel{\partial}{\to} K_1(R, I) \to K_1(R) \to K_1(R/I) \stackrel{\partial'}{\to} K_0(I) \to K_0(R) \to K_0(R/I) \ldots$

Proof: Using snake lemma, see reference [1] for details. ∎

Definition 1.13. ($D_1$ group of a local ring) $D_1$ group of a local ring $A$: $D_1(A)$ is generated by all Dennis-Stein symbols $<a, b>$, where $1 + ab \in A^*; a, b \in A$. Dennis symbol $<a, b>$ is defined by the following with relations:

R1. $<a, b> = -<b, a>$,

R2. $<a, cb> = <ac, b> + <ba, c>$,

R3. $<a, b> + <c, d> = <\theta c \theta^{-1}, \theta d \theta^{-1}> + <a, b>, \theta = (1 + ab)(1 + ba)^{-1}$,

R4. $<a, b> + <\varepsilon a \varepsilon^{-1}, (c + b)\varepsilon^{-1}> = <a, c>, \varepsilon = 1 + ba$.

Definition 1.14. （$D_0$ group of a local ring） $D_0$ group of a local ring $A$: $D_0(A)$ is generated by all Steinberg symbols $\{u,v\}$, where $u,v \in A^*$, which is defined by the following three relations:

V1. $\{u, 1-u\} = 1$,

V2. $\{uv, w\} = \{uvu^{-1}, uwu^{-1}\}\{u, w\}$,

V3. $\{u, v\}\{v, u\} = 1$.

Theorem 1.15. (Magurn) If $F$ is a finite field of characteristic $p$ and $G$ is a finite group whose Sylow $p$-subgroup is a cyclic direct factor, then $K_2(FG) = 0$ if and only if the Sylow $p$-subgroup of $G$ is cyclic.

Proof: See the proof of Theorem1. in reference [2].∎

Definition 1.16. (Hyperelementary Group) A finite group $H$ is hyperelementary if $H$ has a cyclic normal subgroup of prime power index. Denote by $\mathcal{H}(G)$ the set of hyperelementary subgroups of $G$.

Theorem 1.17. (Dress) Let $A$ be a communicative ring, $G$ is a finite group, $H$ is a subgroup of $G$, $f: H \to G$ is a group homomorphism such that we can see $AG$ as a right $AH$-module, which induces a homomorphism between $K$ groups: $f_*: K_n(AH) \to K_n(AG)$. When the morphism runs over every $H \in \mathcal{H}(G)$, we have the isomorphism: $\sum i_* : \varinjlim_{\mathcal{H}(G)} K_2(AH) \to K_2(AG)$, where

$\varinjlim_{\mathcal{H}(G)} K_n(AH) = \oplus_i K_n(AH_i)/\mathcal{R}$, and $\mathcal{R}$ is a subgroup generated by all elements formed in

$(0,\dots,0,x,0,\dots,0,-y,0,\dots 0); f_*(x) = y$, where $f$ is an inclusive map.

Proof: See reference [3].∎

Example 1.18. (Magurn) For alternating group $A_4$, we have: $K_2(F_2[A_4]) \cong Z_2$.

Proof: All the proper subgroups of the alternating group $A_4$ are hyperelementary, and all except $V = \{1, (12)(34), (13)(24), (14)(23)\} \cong Z_2 \oplus Z_2$ are cyclic of order 1, 2 or 3. By Theorem 1.15. $K_2(F_2[Z_2]) = 0, K_2(F_2[Z_3]) = 0$. And in reference[4,Dennis,Keating and Stein], we know that: $K_2(F_2[V]) \cong Z_2^{\oplus 3}$, generated by: $x = \{\sigma, 1+(1+\sigma)\tau\}; y = \{\sigma, \tau\}; z = \{1+(1+\sigma)\tau, \tau\}$. Where $\sigma = (12)(34)$ and $\tau = (13)(24)$. Since there is no inner inclusion from $V$ to any other hyperelementary subgroup of $A_4$, projection to the $K_2(F_2[V])$ summand induces an isomorphism $K_2(F_2[A_4]) \cong \varinjlim_{\mathcal{H}(G)} K_2(F_2[H]) = \oplus_i K_2(F_2[H_i])/\mathcal{R} \cong K_2(F_2[V])/\mathcal{R}$. Where $\mathcal{R}$ is the subgroup generated by all $uv^{-1}$, where $f_*(u) = v$ for some inner inclusion $f: V \to V$. (Since $V \triangleleft A_4$, the group $A_4$ acts on $V$ by conjugation, and $K_2(F_2[V])/\mathcal{R}$ is the module of coinvariants $H_0(A_4, K_2(F_2[V]))$.)

The conjugation action of $A_4$ on $V$ defines a homomorphism $A_4 \to Aut(V)$ with kernel $V$; so the effect of each nontrivial conjugation is the 3-cycle $(\sigma, \tau, \sigma\tau)$ or its inverse. The effect on $K_2(F_2[V])$ takes $x$ to $\{\tau, 1+\sigma\tau+\sigma\} = \{\tau,\tau\}\{\tau,\tau+\sigma+\sigma\tau\} = \{\tau,-\tau\}yz = yz$. By similar computations, we find the effect of this inner inclusion to be the composite of 3-cycles on

$K_2(F_2[V])$: $(x, yz, xyz)(z, xz, xy)$. And the inverse automorphism of $V$ has the inverse effect. So, remembering that $x, y$ and $z$ commute and have order 2, $\mathcal{R} = \{1, xyz, x, yz\}$. Thus $K_2(F_2[V])/\mathcal{R}$ is cyclic of order 2, generated, for instance, by $y = \{\sigma, \tau\}$. So finally we can get: $K_2(F_2[A_4]) \cong K_2(F_2[V])/\mathcal{R} \cong Z_2$. ∎

## 2. Some preliminary studies of the structure of $K_2(F_2[D_4])$

Theorem 2.1. (Tate and Nakayama) Let $p$ is an arbitrary prime number, and $G$ is a $p$-group, then the irreducible representation of $G$ on the fields of characteristic $p$ is only trivial representation.

Proof: See the proof of Theorem2. in Chapter IX of reference [4]. ∎

Corollary 2.2. $F_2[D_4]$ is a local ring.

Proof: Dihedral group $D_4$ is a 2- group, by Theorem 2.1., we know that the irreducible representation of $D_4$ on the field $F_2$ is only trivial representation, and each representation corresponds to an $F_2[D_4]$-module one by one. If $F_2[D_4]$ has two different maxium ideals: $m_1, m_2$, then they can be regarded as $F_2[D_4]$-module. Due to $m_1 \neq m_2$, so there are two different representaions of $D_4$ on the field $F_2$, contradicted! So the proof is completed. ∎

Definition 2.3. (ring with stable range 1) Let $R$ be any ring, if there exists $y \in R$ for every $a, b \in R$ satisify $aR + bR = R$, such that: $a + by \in R^*$, then we call $R$ is a ring with stable range 1. E.g.: $F_2[D_4]$ is a ring with stable range 1.

Definition 2.4. ($p$-vector) Let $R$ be any ring, if $(t_1, t_2, \ldots, t_n) \in R^n$, we can define polynomials with $n$ variables $t_1, t_2, \ldots, t_n$ inductively as follows: $p() = 1; p(t_1) = t_1; p(t_1, t_2) = 1 + t_1 t_2; p(t_1, t_2, t_3) = t_1 + t_3 + t_1 t_2 t_3; \ldots; p(t_1, \ldots, t_n) = p(t_1, \ldots, t_{n-1})t_n + p(t_1, \ldots, t_{n-2})$. When $p(t_1, \ldots, t_n) \in R^*$, we call $(t_1, t_2, \ldots, t_n)$ is a $p$-vector with $n$ variables in $R$, written as: $(t_1, t_2, \ldots, t_n) \in U_n(R)$.

Definition 2.5. (suffix) Let $u = p(t_1, \ldots, t_n) \in R^*$, we call $s = -u^{-1}p(t_1, \ldots, t_{n-1})$ is the suffix of $(t_1, t_2, \ldots, t_n)$. Moreover, let $r = (a_1, \ldots, a_n) \in U_n(R)$, the suffix $s$ of $r$ is also called 1-th suffix of $r$, and the suffix $s'$ of the new vector $(a_2, \ldots, a_n, s)$ is called 2-th suffix of $r$, so we define $n$-th suffix of $r$ inductively. (Remark: Suffixes can be defined here in an inductive way, because every newly generated vector is a $p$- vector, we can see Proposition 3.1.12 in reference [5] for the specific proof.)

Definition 2.6. ($H$-type elements and $W$-type elements) For $\alpha \in R^*$, let: $w_{ij}(\alpha) =$

$x_{ij}(\alpha)x_{ji}(-\alpha^{-1})x_{ij}(\alpha), h_{ij}(\alpha) = w_{ij}(\alpha)w_{ij}(-1)$. For $r = (t_1,\ldots,t_n) \in R^n$, we can define $H$-type elements :

$$H_{ij}(r) = \begin{cases} x_{ij}(t_1)x_{ji}(t_2)\ldots x_{ji}(t_n)x_{ij}(s_1)x_{ji}(s_2), & n \equiv 0 \pmod 2 \\ x_{ij}(t_1)x_{ji}(t_2)\ldots x_{ij}(t_n)x_{ji}(s_1)x_{ij}(s_2)w_{ij}(1), & n \equiv 1 \pmod 2 \end{cases}$$

, where $s_1, s_2$ are 1-th suffixe and 2-th suffixe of $r$ respectively. Then we can define $W$-type elements : $W_{ij}(r) = H_{ij}(r)h_{ij}^{-1}(p(r^*)), r^* = (t_n,\ldots,t_1)$, and let the group $W^{(n)}(R)$ generated by all $W_{ij}(r), r \in U_n(R)$, and let $W(R)$ be the subgroup of $R^*$ generated by all $p(r)p^{-1}(r^*)$. In this way, the mapping is induced: $\phi: W^{(3)}(R) \to W(R): W_{ij}(r) \to p(r)p^{-1}(r^*)$.

**Theorem 2. 7.** There is a short exact sequence: $0 \to K_2(F_2[D_4]) \to W^{(3)}(F_2[D_4]) \xrightarrow{\phi} W(F_2[D_4]) \to 0$.

Proof: Because $F_2[D_4]$ is a ring with stable range 1, so it is immediately obtained from the Theorem 4.1.6 in reference [5].∎

**Theorem 2. 8.** (Kolster) When $R$ is a local ring, we have: $W^{(3)}(R) = W^{(2)}(R) = D_1(R)$, and also have a short exact sequence: $0 \to K_2(R) \to D_1(R) \xrightarrow{\phi} [R^*, R^*] \to 0$.

Proof: See the proof of the main Theorem in reference [6].∎

**Corollary 2. 9.** Group $K_2(F_2[D_4])$ is generated by Dennis-Stein Symbols $<r,s>; r,s \in R, 1 + rs \in R^*, 1 + rt \in R^*, 1 + rst \in R^*, rs = sr, rt = tr, rst = trs, (R = F_2[D_4])$ and they satisfy the following three relations:
DS1: $<r,s> + <s,r> = 0$,
DS2: $<r,s> + <r,t> = <r, s+t+rst>$,
DS3: $<r,st> = <rs,t> + <tr,s>$.

Proof: Because $F_2[D_4]$ is a local ring, so by Theorem 2.8., we know that: $K_2(F_2[D_4]) = Ker\phi = span\{<r,s> | 1 + rs \in (F_2[D_4])^*, p(r,s) = 1 + rs = p(s,r) = 1 + sr\} = span\{<r,s> | 1 + rs \in (F_2[D_4])^*, rs = sr\}$, the proof is completed.∎

**Theorem 2. 10.** $K_2(F_2[D_4])$ is an abelian group and also a 2- group.

Proof: We can see the proof of Theorem5.2.1.(Steinberg) in reference [7] for commutativity. We only need to prove: $K_2(F_2[D_4])$ is a 2- group. We need a lemma:

**Lemma.** Let $A$ be an arbitrary local ring, then there exists a surjective map: $D_0(A) \to D_1(A)$, such that: The kernel of the surjection can be generated by elements that satisfy the following two relations:
V4. $\{u, -u\} = 1$,
V5. $\{u_1, 1 + qu_1\}\{(1 + qu_1)(1 + u_1q)^{-1}u_2(1 + qu_1)^{-1}, (1 + q(u_1 + u_2))(1 + qu_1)^{-1}\}$
$= \{v_1, 1 + qv_1\}\{(1 + qv_1)(1 + v_1q)^{-1}v_2(1 + qv_1)^{-1}, (1 + q(v_1 + v_2))(1 + qv_1)^{-1}\}$.

Where $q \in radA, u_1 + u_2 = v_1 + v_2 \in radA$.

Proof of the lemma: See reference [6] and reference [8] for details.

Back to the original problem, using Corollary 2.9. and the property V2 in the Definition 1.1.14., then we know: $\{r^2, 1+rs\} = \{rrr^{-1}, r(1+rs)r^{-1}\}\{r, 1+rs\} = \{r, 1+rs\}^2$. Therefore, according to the inductive principle, we can get: $\{r, 1+rs\}^n = \{r^n, 1+rs\}$, thus: If $r \in (F_2[D_4])^*$, then we have: $<r,s>^n = \{r, 1+rs\}^n = \{r^n, 1+rs\} \Rightarrow : ord(<r,s>)|ord(r)$. If $r, s \in rad(F_2[D_4])$, due to: $<1,r> = \{1, 1+r\} = \{1, 1+r\}\{r+1, r+1\} = \{r+1, r+1\} = 1 = <r,1>$, then use the DS2 in Corollary 2.9. $\Rightarrow : <r,s> = <r,s><r,1> = <r, 1+s+rs> = <1+s+rs, r>^{-1}$, at this time: $1+s+rs \in (F_2[D_4])^*$, as all mentioned above, we know: $K_2(F_2[D_4])$ is generated by all $<r,s>$ which satisfy all conditions in the Lemma above withing: $r \in (F_2[D_4])^*$. Due to $ord((F_2[D_4])^*) = 128 = 2^7$, so $ord(r)$ is also a power of 2, and notice: $ord(<r,s>)|ord(r)$, then we know: $ord(<r,s>)$ is a power of 2, so, $K_2(F_2[D_4])$ Is a 2- group, the proof is completed.■

Theorem 2. 11. The minimum direct summand terms of $K_2(F_2[D_4])$ can only be $Z_2$ and $Z_4$.

Proof: Notice that there is a short exact sequence: $0 \to (\sigma^2 - 1)F_2[D_4] \to F_2[D_4] \to F_2[Z_2 \oplus Z_2] \to 0$, and $(\sigma^2 - 1)F_2[D_4]$ is nilpotent $((\sigma^2 - 1) \subseteq Cent(F_2[D_4]), (\sigma^2 - 1)^2 = 0)$, and the order of invertible elements in the local ring $F_2[Z_2 \oplus Z_2]$ can only be 1 or 2, so we know the order of invertible elements in $F_2[D_4]$ can only be 1,2,4. The proof is completed.■

# 3. Construct the structure of the group $K_2(F_2[D_4])$

Theorem 3. 1. $K_2(F_2[D_4]) \cong Z_2$。

Proof: Let the ring $R = F_2[D_4]$, and consider a surjective map between groups: $D_4 = <\sigma, \tau | \sigma^4 = \tau^2 = 1, \sigma\tau = \tau\sigma^3> \to Z_2 \oplus Z_2 = <\sigma_1, \tau | \sigma_1^2 = \tau^2 = 1, \sigma_1\tau = \tau\sigma_1> : \sigma \to \sigma_1, \tau \to \tau$, Which induces a surjective map between two group algebras: $f: R \to F_2[Z_2 \oplus Z_2] \cong R/I$, where $I$ is the kernel of $f$, and also is the center of $R$, i.e.: $I = (\sigma^2 - 1)R$. Therefore, we have a short exact sequence: $0 \to I \to R \to R/I \to 0$, then let the functor $K_2$ acts on this short exact sequence, and due to Theorem6.2 in reference [1] and Definition 1.1.10. (the long exact sequence theorem of K-theory) and combine with Definition 1.11. (In fact, all the relative $K_2$ groups in this paper are just the kernel of the map of the relative $K_2$ group in reference [1] to the $K_2$ group), so we can get the long exact sequence immediately: $0 \to K_2(R,I) \to K_2(R) \xrightarrow{f} K_2(R/I) \xrightarrow{\partial} K_1(R,I) \xrightarrow{g} K_1(R) \xrightarrow{h} K_1(R/I) \xrightarrow{\partial'} K_0(R,I) \to K_0(R) \to K_0(R/I) \to \ldots$.     Then,let's consider the first 4 items in the long exact sequence: $0 \to K_2(R,I) \to K_2(R) \xrightarrow{f} K_2(R/I)$. We define $<R,R>$ is a group generated by all Denis-Stein symbols which are: $<r,s> ; r, s \in R; 1+rs \in R^*, 1+rt \in R^*, 1+rst \in R^*; rs = sr, rt = tr, rst = trs$, and the following three relations are satisfied at the same time:

DS1: $<r,s>+<s,r>=0$,
DS2: $<r,s>+<r,t>=<r,s+t+rst>$,
DS3: $<r,st>=<rs,t>+<tr,s>$.

By Corollary 2.9. we know that: $K_2(R) \cong <R,R>$, notice that: $I^2 = 0$, so: For arbitrary $r_1, r_2 \in R$, satisfy: $r_1 r_2 = r_2 r_1, 1 + r_1 r_2 \in R^*$, and for any $i \in I$, we have: $<r_1, r_2> + <r_1, i(1+r_1 r_2)^{-1}> = <r_1, r_2 + i(1+r_1 r_2)^{-1} + r_1 r_2 i(1+r_1 r_2)^{-1}> = <r_1, r_2 + (1+r_1 r_2)i(1+r_1 r_2)^{-1}> = <r_1, r_2 + i>$, combining with the definition of the relative $K_2$ group (Definition 1.10.), we have $K_2(R)/<<i,r> | i \in I, r \in R, <i,r> \in <R,R>>> \subseteq <R/I, R/I> \cong K_2(R/I) \cong K_2(R)/K_2(R,I)$, so: $K_2(R,I) \subseteq <<i,r> | i \in I, r \in R, <i,r> \in <R,R>>>$. Then, we will prove a lemma:

**Lemma 1**  $<<i,j> | i, j \in I, <i,j> \in <R,R>>> = 0$。

Proof of lemma 1: For $\forall i, j \in I$, we let $i = (\sigma^2 - 1)a, j = (\sigma^2 - 1)b; a, b \in R$, by DS3, we have: $<i,j> = <(\sigma^2-1)a, (\sigma^2-1)b> = <\sigma^2-1, a(\sigma^2-1)b> - <(\sigma^2-1)b(\sigma^2-1), a> = <\sigma^2-1, (\sigma^2-1)ab> = <(\sigma-1)^2, (\sigma-1)^2 ab>$, notice that: $I^2 = 0$, so: $\forall r = a_1 + a_2\sigma + a_3\sigma^2 + a_4\sigma^3 + b_1\tau + b_2\sigma\tau + b_3\sigma^2\tau + b_4\sigma^3\tau \in R; a_1 \sim a_4, b_1 \sim b_4 \in \{0,1\}$, we all have that: $(\sigma^2-1)r = (\sigma^2-1)(a_1 + a_2\sigma + a_3\sigma^2 + a_4\sigma^3 + b_1\tau + b_2\sigma\tau + b_3\sigma^2\tau + b_4\sigma^3\tau) = (\sigma^2-1)(a_1 + a_2\sigma + b_1\tau + b_2\sigma\tau) = (\sigma^2-1)(a_1 + a_2\sigma + (b_1 + b_2\sigma)\tau)$ ($*$), so we can get a conclusion: $<i,j> = <(\sigma^2-1)a, (\sigma^2-1)b> = <(\sigma^2-1)(a_1 + a_2\sigma + (b_1 + b_2\sigma)\tau), (\sigma^2-1)(a'_1 + a'_2\sigma + (b'_1 + b'_2\sigma)\tau)> = <(\sigma^2-1)(a_1 + a_2\sigma + (b_1 + b_2\sigma)\tau)(\sigma-1), (\sigma-1)(a'_1 + a'_2\sigma + (b'_1 + b'_2\sigma)\tau)> + <(\sigma-1)^3(a'_1 + a'_2\sigma + (b'_1 + b'_2\sigma)\tau)(a_1 + a_2\sigma + (b_1 + b_2\sigma)\tau), \sigma-1> = <(\sigma-1)^3(a'_1 + a'_2\sigma + (b'_1 + b'_2\sigma)\tau)(a_1 + a_2\sigma + (b_1 + b_2\sigma)\tau), \sigma-1> + <(\sigma-1)^3(a'_1 + a'_2\sigma + (b'_1 + b'_2\sigma)\tau)(a_1 + a_2\sigma + (b_1 + b_2\sigma)\tau), \sigma-1> = 2<(\sigma-1)^3(a'_1 + a'_2\sigma + (b'_1 + b'_2\sigma)\tau)(a_1 + a_2\sigma + (b_1 + b_2\sigma)\tau), \sigma-1> = 2<(\sigma-1)^3(k_1 + k_2\sigma + (l_1 + l_2\sigma)\tau), \sigma-1> = 2<\sigma-1, (\sigma-1)^3(k_1 + k_2\sigma + (l_1 + l_2\sigma)\tau)> = 2<\sigma-1, (\sigma-1)^3(k_1 + k_2\sigma)> + 2<\sigma-1, (\sigma-1)^3(l_1 + l_2\sigma)\tau>$ (#) The reason why the antepenultimate equal sign is established in here is due to the formula ($*$), where $k_1, k_2, l_1, l_2 \in \{0,1\}$. And, we notice that: If $k_1 = k_2 = 1$, we immediately have: $<\sigma-1, (\sigma-1)^3(k_1 + k_2\sigma)> = 0$; if $k_1 = 0, k_2 = 1$, then have: $<\sigma-1, (\sigma-1)^3(k_1 + k_2\sigma)> = <\sigma-1, (\sigma^2-\sigma)(\sigma^2-1)> = <\sigma-1, (\sigma^2-1-\sigma+1)(\sigma^2-1)> = <\sigma-1, (\sigma-1)^3>$; if $k_1 = 1, k_2 = 0$, then have: $<\sigma-1, (\sigma-1)^3(k_1 + k_2\sigma)> = <\sigma-1, (\sigma-1)^3>$. With similarly: If $l_1 = l_2 = 1$, immediately we have: $<\sigma-1, (\sigma-1)^3(l_1 + l_2\sigma)\tau> = 0$; if $l_1 = 0, l_2 = 1$, we have: $<\sigma-1, (\sigma-1)^3(l_1 + l_2\sigma)\tau> = <\sigma-1, (\sigma-1)^3\sigma\tau> = <\sigma-1, (\sigma^2-\sigma)(\sigma^2-1)\tau> = <\sigma-1, (\sigma^2-1-\sigma+1)(\sigma^2-1)\tau> = <\sigma-1, (\sigma-1)^3\tau>$; if $l_1 = 1, l_2 = 0$, then have: $<\sigma-1, (\sigma-1)^3(l_1 + l_2\sigma)\tau> = <\sigma-1, (\sigma-1)^3\tau>$. So, substituting the above results into the formula (#), we immediately get: $<<i,j> | i, j \in I, <i,j> \in <R,R>>> = \{0, 2<\sigma-1, (\sigma-1)^3>, 2<\sigma-1, (\sigma-1)^3\tau>, 2<\sigma-1, (\sigma-1)^3> + 2<\sigma-1, (\sigma-1)^3\tau>\}$. Using DS2, we find that: $2<\sigma-1, (\sigma-1)^3> = <\sigma-1, (\sigma-1)^7> = <\sigma-1, 0> = 0$ and $2<\sigma-1, (\sigma-1)^3\tau> = <\sigma-1, (\sigma-1)^4\tau(\sigma-1)^3\tau> = <\sigma-1, 0> = 0$, so as all mentioned above, we get: $<<i,j> | i, j \in I, <i,j> \in <R,R>>> = 0$, the proof of this lemma is completed.

Back to the original problem, due to $K_2(R,I) \subseteq <<i,r> | i \in I, r \in R, <i,r> \in <R,R>>>$ and $I = Cent(R) = Ker(f)$, we know: $K_2(R,I)$ is generated by all $<x,y> \in <R,R>$, where $x \in I$ or $y \in I$, and due to the formula ($*$) in lemma 1, we can get that: $\{<(\sigma^2-1)(a_0 + a_1(\sigma-1) + a_2(\tau-1) + a_3(\sigma-1)(\tau-1), b_0 + b_1(\sigma-1) + b_2(\tau-1) + b_3(\sigma-1)(\tau-1)>$

$|\forall a_i, b_i \in \{0,1\}\}$ is a set of generators of the group $K_2(R,I)$, according to the conclusion of lemma 1.: $<<i,j>|i,j \in I, <i,j>\in<R,R>>= 0$, we can get: $<i,x+y>=<i,x+y>+<i,ixy>=<i,x+y+ixy>=<i,x>+<i,y>$, $i \in I; x,y \in R$, thus we have: The set of all generators of $K_2(R,I)$ is: $\{<(\sigma^2-1)(\sigma-1)^i(\sigma-1)^j, (\sigma-1)^{i'}(\tau-1)^{j'}>|i,j,i',j' \in \{0,1\}\}$, there are $2^4 = 16$ elements in this set, combining with the conclusion of lemma 1. that: $<<i,j>|i,j \in I, <i,j>\in<R,R>>= 0$, the set can be simplified to the following 4 elements: $<(\sigma-1)^2, \tau-1>$, $<(\sigma-1)^2, \sigma-1>$, $<(\sigma-1)^2(\tau-1), \sigma-1>$, $<(\sigma-1)^3, \tau-1>$; then we immediately find that:

$<(\sigma-1)^3, \tau-1>=<\tau-1, (\sigma-1)^3>=<(\tau-1)(\sigma-1)^2, \sigma-1>+<(\sigma-1)(\tau-1), (\sigma-1)^2>=<(\sigma-1)^2(\tau-1), \sigma-1>+<(\sigma-1)^2, (\sigma-1)(\tau-1)>=<(\sigma-1)^2(\tau-1), \sigma-1>+0=<(\sigma-1)^2(\tau-1), \sigma-1>$, and: $<(\sigma-1)^2, \sigma-1>=<\sigma-1, (\sigma-1)^2>=<(\sigma-1)^2, \sigma-1>+<(\sigma-1)^2, \sigma-1>=0$. Therefore, as all mentioned above, we know that: $K_2(R,I)$ can be generated by $<(\sigma-1)^2, \tau-1>$ and $<(\sigma-1)^3, \tau-1>$. Combining with the conclusion of lemma 1 that: $<<i,j>|i,j \in I, <i,j>\in<R,R>>= 0$, we know the generators of $K_2(R)$, can be generated by the generators of the group $K_2(R,I)$: $<(\sigma-1)^2, \tau-1>$, $<(\sigma-1)^3, \tau-1>$ and by a set of generators of the group $S = <<x,y>|x,y \in R \setminus I, <x,y>\in<R,R>>$. Next, we will determine a group of generators of the group $S = <<x,y>|x,y \in R \setminus I, <x,y>\in<R,R>>$: Known from the proof process of lemma 1: The elements on both sides of the generators of the group $S$ are as follows: $a_0 + a_1\sigma + a_2\tau + a_3\sigma\tau$, where $a_0 \sim a_3 \in F_2$. Now let: $x = a_0 + a_1\sigma + a_2\tau + a_3\sigma\tau = a_0 + x_1, y = b_0 + b_1\sigma + b_2\tau + b_3\sigma\tau = b_0 + y_1$, due to: $xy = (a_0+x_1)(b_0+y_1) = a_0b_0 + a_0y_1 + b_0x_1 + x_1y_1, yx = (b_0+y_1)(a_0+x_1) = b_0a_0 + b_0x_1 + a_0y_1 + y_1x_1$, so: $xy = yx \Leftrightarrow x_1y_1 = y_1x_1$, moreover we have: $x_1y_1 = (a_1\sigma + a_2\tau + a_3\sigma\tau)(b_1\sigma + b_2\tau + b_3\sigma\tau) = a_1b_1\sigma^2 + a_1b_2\sigma\tau + a_1b_3\sigma^2\tau + a_2b_1\sigma^3\tau + a_2b_2 + a_2b_3\sigma^3 + a_3b_1\tau + a_3b_2\sigma + a_3b_3, y_1x_1 = (b_1\sigma + b_2\tau + b_3\sigma\tau)(a_1\sigma + a_2\tau + a_3\sigma\tau) = a_2b_2 + a_3b_3 + a_3b_2\sigma + a_1b_1\sigma^2 + a_2b_3\sigma^3 + a_3b_1\tau + a_1b_2\sigma\tau + a_1b_3\sigma^2\tau + a_2b_1\sigma^3\tau$, Through comparing coefficients, we get: $x_1y_1 = y_1x_1 \Leftrightarrow a_1b_2 = a_2b_1, a_3b_1 = a_1b_3, a_3b_2 = a_2b_3$, because $x_1, y_1$ are in the set: $\{\sigma, \tau, \sigma\tau, \sigma+\tau, \sigma+\sigma\tau, \tau+\sigma\tau, \sigma+\tau+\sigma\tau\}$, satisfying both the communicativity and $1+x_1y_1 \in R^*$ are only the following 3 generators: $<\sigma+\tau, \sigma+\tau>$, $<\sigma+\sigma\tau, \sigma+\sigma\tau>$, $<\tau+\sigma\tau, \tau+\sigma\tau>$, so as mentioned all above: The following 8 elements: $<(\sigma-1)^2, \tau-1>$, $<(\sigma-1)^3, \tau-1>$, $<\sigma+\tau, \sigma+\tau>$, $<\sigma+\sigma\tau, \sigma+\sigma\tau>$, $<\tau+\sigma\tau, \tau+\sigma\tau>$, $<\sigma+\tau, 1+\sigma+\tau>$, $<\sigma+\sigma\tau, 1+\sigma+\sigma\tau>$, $<\tau+\sigma\tau, 1+\tau+\sigma\tau>$, can be a group of generators of the group $K_2(R)$. We notice that for any arbitrary $x \in \{\sigma+\tau, \sigma+\sigma\tau, \tau+\sigma\tau\}$, by DS2, we can get: $<x, 1+x+x^2>=<x,1>+<x,x>=<x,x>$, by DS3 and DS1 we can get: $-<x^2,x>=<x,x^2>=<x^2,x>+<x^2,x> \Rightarrow : 3<x^2,x>=0 \Rightarrow : <x,x^2>=0$, by DS2 again and by $<x,x^2>=0$, we can get: $<x,1+x>=<x,1+x>+<x,x^2>=<x,1+x+x^2+x(1+x)x^2>=<x,1+x+x^2+x^3>=<x,1+x+x^2>+<x,x^3>$, and by DS2 we can also get: $<x,x^3>=<x,x+x+x^3>=<x,x>+<x,x>=0$, so as mentioned all above, we have: $<x,1+x>=<x,1+x+x^2>=<x,x>$, i.e.: A set of generators of $K_2(R)$ can be simplified to the following 5 elements: $<(\sigma-1)^2, \tau-1>$, $<(\sigma-1)^3, \tau-1>$, $<\sigma+\tau, \sigma+\tau>$, $<\sigma+\sigma\tau, \sigma+\sigma\tau>$, $<\tau+\sigma\tau, \tau+\sigma\tau>$. In order to simplify the number of generators further, we consider the map in the long exact sequence of K-theory: $K_2(R) \xrightarrow{f} K_2(R/I)$, it is known from reference [2] that: $K_2(R/I)$ is generated by the following 3 generators: $\alpha = \{\sigma^2, 1+(1+\sigma^2)\tau\} = <\sigma^2, \tau+\sigma^2\tau>$, $\beta = \{\sigma^2, \tau\} = <\sigma^2, \sigma^2\tau+\sigma^2>$, $\gamma = \{1+(1+\sigma^2)\tau, \tau\} = <$

$\sigma^2\tau + \tau + 1, \sigma^2\tau + \sigma^2 >$, notice that: $f(<\sigma+\tau,\sigma+\tau>) = <\sigma^2+\tau,\sigma^2+\tau>$, and $<\sigma^2+\tau,\sigma^2+\tau> + <\sigma^2+\tau,\sigma^2> = <\sigma^2+\tau,\sigma^2+\tau+\sigma^2+(\sigma^2+\tau)(\sigma^2+\tau)\sigma^2> = <\sigma^2+\tau,\tau>$, due to that: $<\sigma^2+\tau,\sigma^2> + <\sigma^2+\tau,\tau> = <\sigma^2+\tau,\sigma^2+\tau+(\sigma^2+\tau)\sigma^2\tau> = <\sigma^2+\tau,0> = 0$, and $ord(<\sigma^2+\tau,\sigma^2>) = ord(<\sigma^2+\tau,\tau>) = 2$, so we get: $<\sigma^2+\tau,\sigma^2> = <\sigma^2+\tau,\tau>$, from this, we can get: $f(<\sigma+\tau,\sigma+\tau>) = <\sigma^2+\tau,\sigma^2+\tau> = 0$. Similarly: $f(<\sigma+\sigma\tau,\sigma+\sigma\tau>) = <\sigma^2+\sigma^2\tau,\sigma^2+\sigma^2\tau> = <\sigma^2(\tau+1),\sigma^2(\tau+1)> = <\sigma^2(\tau+1)\sigma^2,\tau+1> + <(\tau+1)\sigma^2(\tau+1),\sigma^2> = <\tau+1,\tau+1> + <0,\sigma^2> = <\tau+1,\tau+1>$, due to: $<\tau+1,\tau+1> = <\tau+1,\tau+1> + <\tau+1,1> + <\tau+1,1> = <\tau+1,\tau> + <\tau+1,1> = <\tau+1,\tau+1+\tau+1> = <\tau+1,0> = 0$, so we get: $f(<\sigma+\sigma\tau,\sigma+\sigma\tau>) = <\sigma^2+\sigma^2\tau,\sigma^2+\sigma^2\tau> = 0$, and: $f(<\tau+\sigma\tau,\tau+\sigma\tau>) = <\tau+\sigma^2\tau,\tau+\sigma^2\tau> = <\tau(\sigma^2+1),\tau(\sigma^2+1)> = <\sigma^2+1,\sigma^2+1> = 0$. As all mentioned above: A group of generators of the $K_2(R)$ can be further simplified as: $<(\sigma-1)^2,\tau-1>$, $<(\sigma-1)^3,\tau-1>$, i.e.: A group of generators of the $K_2(R,I)$, so now we get: $K_2(R) \cong K_2(R,I)$, there are no elements of order 4 in the $K_2(R)$, and there are only 3 possible structures of the $K_2(R)$: $0, Z_2, Z_2 \oplus Z_2$. The following two short exact sequences are known from Thm3.7(Kolster,1984) and Thm3.9(Kolster,1985) in reference [9]: $0 \to K'_2(2,R) \to D_1(R) \to E_2(R) \cap GL_1(R) \to 0$ and $0 \to K_2(R) \to D_1(R) \to [R^*,R^*] \to 0$, because $R$ is a local ring, so due to reference [10, M.Kolster,1982], we know that: $K_2(R) \cong K'_2(2,R)$, where the definition of $K'_2(2,R)$ can be seen in Def.1.6 in reference [9], $D_1(R)$ is a group generated by all Dennis-Stein symbols given in the Part 1, that: $<a,b>$; $a,b \in R$, $1+ab \in R^*$ (Remark: $a,b$ are needn't to be communicated), and the generators satisfie the 4 relations in Definition 1.13.:

R1. $<a,b> = -<b,a>$,

R2. $<a,cb> = <ac,b> + <ba,c>$,

R3. $<a,b> + <c,d> = <\theta c\theta^{-1},\theta d\theta^{-1}> + <a,b>$, $\theta = (1+ab)(1+ba)^{-1}$,

R4. $<a,b> + <\varepsilon a\varepsilon^{-1},(c+b)\varepsilon^{-1}> = <a,c>$, $\varepsilon = 1+ba$.

Due to reference [10, M.Kolster,On injective stability for K2, 1982], we know that: $E_2(R) \cap GL_1(R) \cong [R^*,R^*]$, now we will construct the elements in $E_2(R) \cap GL_1(R)$:

Let $a,b \in R$, $1+ab \in R^*$, consider the operations in the elementary matrices group $E_2(R)$, that: $\begin{pmatrix} 1 & a \\ 0 & 1 \end{pmatrix}\begin{pmatrix} 1 & 0 \\ b & 1 \end{pmatrix} = \begin{pmatrix} 1+ab & a \\ b & 1 \end{pmatrix}$, $\begin{pmatrix} 1 & 0 \\ b(1+ab)^{-1} & 1 \end{pmatrix}\begin{pmatrix} 1+ab & a \\ b & 1 \end{pmatrix}\begin{pmatrix} 1 & (1+ab)^{-1}a \\ 0 & 1 \end{pmatrix} = \begin{pmatrix} 1+ab & 0 \\ 0 & 1+b(1+ab)^{-1}a \end{pmatrix}$, so we can get that: $\begin{pmatrix} 1+ab & 0 \\ 0 & 1+b(1+ab)^{-1}a \end{pmatrix} \in E_2(R)$, and we notice there are some facts: $\begin{pmatrix} (1+ab)^3 & ab \\ 1+ab+(ab)^2 & 1 \end{pmatrix} = \begin{pmatrix} 1 & ab \\ 0 & 1 \end{pmatrix}\begin{pmatrix} 1 & 0 \\ 1+ab+(ab)^2 & 1 \end{pmatrix} \in E_2(R)$, $\begin{pmatrix} (1+ab)^3 & ab \\ 1+ab+(ab)^2 & 1 \end{pmatrix}\begin{pmatrix} 1+ab & 0 \\ 0 & 1+b(1+ab)^{-1}a \end{pmatrix} = \begin{pmatrix} 1 & ab+ab^2(1+ab)^{-1}a \\ 1+(ab)^3 & 1+b(1+ab)^{-1}a \end{pmatrix} = \begin{pmatrix} 1 & u \\ v & 1+b(1+ab)^{-1}a \end{pmatrix} \in E_2(R)$.

So, we can get: $\begin{pmatrix} 1 & 0 \\ v(1+uv)^{-1} & 1 \end{pmatrix}\begin{pmatrix} 1 & u \\ v & 1+b(1+ab)^{-1}a \end{pmatrix}\begin{pmatrix} 1 & (1+uv)^{-1}u \\ 0 & 1 \end{pmatrix} = \begin{pmatrix} 1 & 0 \\ 0 & (1+ab)(1+ba)^3 \end{pmatrix} \in E_2(R) \cap GL_1(R)$. Then, we will prove a lemma:

Lemma 2 $[R^*,R^*] = \{1, \sigma^2, \sigma^2+\tau+\sigma\tau+\sigma^2\tau+\sigma^3\tau, 1+\tau+\sigma\tau+\sigma^2\tau+\sigma^3\tau, \sigma+\sigma^2+\sigma^3+$

$\sigma\tau + \sigma^3\tau, 1 + \sigma + \sigma^3 + \sigma\tau + \sigma^3\tau, 1 + \sigma + \sigma^3 + \tau + \sigma^2\tau, \sigma + \sigma^2 + \sigma^3 + \tau + \sigma^2\tau\} \cong Z_2 \oplus Z_2 \oplus Z_2$.

Proof of lemma 2: We notice that: $(\sigma - 1)^2 R + 1 \subseteq Cent(R^*), |(\sigma - 1)^2 R + 1| = 16$, so we have : $|[R^*, R^*]| = |E_2(R) \cap GL_1(R)| \le |(\sigma - 1)^2 R + 1| = 16$, and we also find that: $a = \sigma + 1, b = \sigma, (1 + ab)(1 + ba)^3 = \sigma + \sigma^2 + \sigma^3 + \sigma\tau + \sigma^3\tau; a = \sigma + \tau, b = \sigma, (1 + ab)(1 + ba)^3 = \sigma^2; a = \sigma^2 + \tau, b = \sigma\tau, (1 + ab)(1 + ba)^3 = \sigma + \sigma^2 + \sigma^3 + \tau + \sigma^2\tau; a = \sigma^2 + \sigma\tau, b = \sigma\tau, (1 + ab)(1 + ba)^3 = \sigma^2 + \tau + \sigma\tau + \sigma^2\tau + \sigma^3\tau; a = \sigma + \tau, b = \sigma\tau, (1 + ab)(1 + ba)^3 = 1 + \sigma + \sigma^3 + \tau + \sigma^2\tau; a = \sigma + \sigma\tau, b = \sigma^2\tau, (1 + ab)(1 + ba)^3 = 1 + \sigma + \sigma^3 + \sigma\tau + \sigma^3\tau$ , the six pairwise different $(1 + ab)(1 + ba)^3$ above, the values of them are all in the group: $E_2(R) \cap GL_1(R) \cong [R^*, R^*]$, so we get that: $(1 + \sigma + \sigma^3 + \tau + \sigma^2\tau)(1 + \sigma + \sigma^3 + \sigma\tau + \sigma^3\tau) = 1 + \tau + \sigma\tau + \sigma^2\tau + \sigma^3\tau \in [R^*, R^*]$, and $1 \in [R^*, R^*]$, now there are 8 elements all above. Now, we only need to show: $1 + \tau + \sigma^2\tau \in \{(\sigma - 1)^2 R + 1\} \setminus [R^*, R^*]$. If it is not correct, then: $\exists a, b \in R^*, s.t.: aba^{-1}b^{-1} = 1 + \tau + \sigma^2\tau \Leftrightarrow ab + ba = (\sigma^2\tau + \tau)ba$ , due to the fact: $(\sigma - 1)^2 R + 1 \subseteq Cent(R^*), |(\sigma - 1)^2 R + 1| = 16, |R^*| = 128$ , so we have: $[R^*, R^*] = [1 + K, 1 + K]; K = k_1(\sigma - 1) + k_2(\tau - 1) + k_3(\sigma\tau - 1), k_1, k_2, k_3 \in F_2$, so we can assume：$a = 1 + a_1(\sigma - 1) + a_2(\tau - 1) + a_3(\sigma\tau - 1), b = 1 + b_1(\sigma - 1) + b_2(\tau - 1) + b_3(\sigma\tau - 1), a_1, a_2, a_3, b_1, b_2, b_3 \in F_2$, now, we have that: $ab + ba = (a_1 b_2 + a_2 b_1 + b_2 a_3 + b_3 a_2) + (a_1 b_2 + b_1 a_2 + b_3 a_2)\sigma + b_2 a_3 \sigma^3 + (a_1 b_2 + b_1 a_2 + b_2 a_3 + a_3 b_1 + b_3 a_1 + b_3 a_2)\tau + (a_1 b_2 + b_2 a_3 + b_3 a_2)\sigma\tau + (a_1 b_3 + b_1 a_3)\sigma^2 \tau + a_2 b_1 \sigma^3 \tau$, and we also have a fact is that: $(\sigma^2 \tau + \tau)ba = (a_2 + b_1 a_2 + b_1 a_3 + b_2 + b_2 a_1 + b_3 a_1 + b_3 a_2) + (a_3 + b_1 a_2 + b_1 a_3 + b_2 a_1 + b_3 a_1 + b_3 a_2)\sigma + (a_2 + b_1 a_2 + b_1 a_3 + b_2 + b_2 a_1 + b_3 a_1 + b_3 a_2)\sigma^2 + (a_3 + b_1 a_2 + b_1 a_3 + b_2 a_1 + b_3 a_1 + b_3 a_2)\sigma^3 + (1 + a_1 + a_2 + a_3 + b_1 + b_1 a_2 + b_1 a_3 + b_2 + b_2 a_1 + b_3 a_1 + b_3 a_2)\tau + (a_1 + b_1 + b_1 a_2 + b_1 a_3 + b_2 a_1 + b_3 a_1 + b_3 a_2)\sigma\tau + (1 + a_1 + a_2 + a_3 + b_1 + b_1 a_2 + b_1 a_3 + b_2 + b_2 a_1 + b_3 a_1 + b_3 a_2)\sigma^2 \tau + (a_1 + b_1 + b_1 a_2 + b_1 a_3 + b_2 a_1 + b_3 a_1 + b_3 a_2)\sigma^3 \tau$ , by comparing the coefficients, we get 8 homogeneous linear equations as follow:

$$\begin{cases} a_2 + b_1 a_3 + b_2 + b_3 a_1 + b_2 a_3 = 0 & (1) \\ a_3 + b_1 a_3 + b_3 a_1 = 0 & (2) \\ a_2 + b_1 a_2 + b_1 a_3 + b_2 + b_2 a_1 + b_3 a_1 + b_3 a_2 = 0 & (3) \\ a_3 + b_1 a_2 + b_1 a_3 + b_2 a_1 + b_3 a_1 + b_3 a_2 + b_2 a_3 = 0 & (4) \\ 1 + a_1 + a_2 + a_3 + b_1 + b_2 + b_2 a_3 = 0 & (5) \\ a_1 + b_1 + b_1 a_2 + b_1 a_3 + b_3 a_1 + b_2 a_3 = 0 & (6) \\ 1 + a_1 + a_2 + a_3 + b_1 + b_1 a_2 + b_2 + b_2 a_1 + b_3 a_2 = 0 & (7) \\ a_1 + b_1 + b_1 a_3 + b_2 a_1 + b_3 a_1 + b_3 a_2 = 0 & (8) \end{cases}$$

By substituting (2) into (5), (6), we can get: $a_1 + b_1 + b_1 a_2 + a_3 + b_2 a_3 = 0$ ($*$), and by substituting ($*$) into (5), we get: $1 + a_2 + b_2 + b_1 a_2 = 0$ ($**$), by substituting (1) into (3), we can get: $b_1 a_2 + b_2 a_1 + b_3 a_2 + b_2 a_3 = 0$ ($***$). If: $a_2 = 1$, then due to ($**$), we immediately get: $b_1 = b_2$, combining with (5), we also can get: $a_1 + a_3 + b_2 a_3 = 0 = a_1 + a_3 + b_1 a_3$, thus, we have: $a_1 = (1 + b_1) a_3$; (i): $1 = a_1 = (1 + b_1)a_3 \Rightarrow : b_1 = 0, a_3 = 1$, by substituting it into ($*$), we get: $a_1 + 1 + b_2 = 0 \Rightarrow : a_1 + b_2 = 1 = a_1 + b_1$, substitute the result just obtained into ($*$), we get: $b_1 = b_2 = 0, a_1 = a_2 = 1$, now, substitute the result just obtained into (8), we get: $1 = 0$, contradicted! (ii): $0 = a_1 = (1 + b_1)a_3$, substitute it into (7), we get: $a_3 + b_1 + b_2 + b_3 = 0$, substitute $0 = a_1 = (1 + b_1)a_3$ into (2), we obtain: $a_3 = b_1 a_3$, substitute $0 = a_1 = (1 + b_1)a_3$ into (3), we obtain: $1 + b_1 + a_3 + b_3 = 0$, so we conclude that: $b_2 = 1$, then substitute $0 = a_1 = (1 + b_1)a_3$ and $b_2 = 1$ into (8), we obtain: $b_1 + a_3 + b_3 = 0$, it contradicts to $1 + b_1 + a_3 + b_3 = 0$! So, it can only be: $a_2 = 0$, substitute it into ($**$) we get: $b_2 = 1$, by substituting $a_2 = 0, b_2 = 1$ into $(*)$ and $(***)$, we get: $a_1 + a_3 = 0 = b_1 + a_3 \Rightarrow : b_1 = a_3$, substitute $a_2 = 0, b_2 = 1, b_1 = a_3$ into (8), we get: $b_3 a_1 = 0$, then substitute $a_2 = 0, b_2 = 1, b_1 = a_3, a_1 + a_3 = 0 = b_1 + a_3, b_3 a_1 = 0$ into (6), we get: $b_1 a_3 = a_3$, at last, by substituting all the results mentioned above into (3), we obtain that: $a_3 + 1 + a_1 = 1 = 0$, also contradicted! So, as mentioned all above, we have: $E_2(R) \cap GL_1(R) \cong [R^*, R^*] = \{1, \sigma^2, \sigma^2 + \tau + \sigma\tau + \sigma^2\tau + \sigma^3\tau, 1 + \tau + \sigma\tau + \sigma^2\tau + \sigma^3\tau, \sigma + \sigma^2 + \sigma^3 + \sigma\tau + \sigma^3\tau, 1 + \sigma + \sigma^3 + \sigma\tau + \sigma^3\tau, 1 + \sigma + \sigma^3 + \tau + \sigma^2\tau, \sigma + \sigma^2 + \sigma^3 + \tau + \sigma^2\tau\} \cong Z_2 \oplus Z_2 \oplus Z_2$. The proof of Lemma 2 is completed.

Back to the original problem, and due to the conclusion of lemma 2: $E_2(R) \cap GL_1(R) \cong [R^*, R^*] = \{1, \sigma^2, \sigma^2 + \tau + \sigma\tau + \sigma^2\tau + \sigma^3\tau, 1 + \tau + \sigma\tau + \sigma^2\tau + \sigma^3\tau, \sigma + \sigma^2 + \sigma^3 + \sigma\tau + \sigma^3\tau, 1 + \sigma + \sigma^3 + \sigma\tau + \sigma^3\tau, 1 + \sigma + \sigma^3 + \tau + \sigma^2\tau, \sigma + \sigma^2 + \sigma^3 + \tau + \sigma^2\tau\} \cong Z_2 \oplus Z_2 \oplus Z_2 \subseteq (\sigma^2 - 1)R + 1 = Cent(R^*)$, then we immediately have: $\theta = (1 + ab)(1 + ba)^{-1} \in Cent(R^*)$, thus R3. can be rewritten as: $<a,b> + <c,d> = <c,d> + <a,b>$, i.e.: $D_1(R)$ is an abelian group, and we have a short exact sequence: $0 \to K_2(R) \to D_1(R) \to Z_2^{\oplus 3} \to 0$. Because $K_2(R)$ is a subgroup of $D_1(R)$, then the following operations are all meaningful: $<(\sigma-1)^2, \tau-1> = <(\sigma-1)^2, \tau> = <\tau, (\sigma-1)^2> = <\tau(\sigma-1), \sigma-1> + <(\sigma-1)\tau, \sigma-1>$, it concludes that: $<(\sigma-1)^2, \tau-1> = -<(\sigma-1)^2, \tau-1> = -<\tau(\sigma-1), \sigma-1> - <(\sigma-1)\tau, \sigma-1> = <\sigma-1, (\sigma-1)\tau> + <\sigma-1, \tau(\sigma-1)>$, using the property R4. and Lemma 1, we can get the following formula:

$<\sigma-1, (\sigma-1)\tau> - <\sigma-1, \tau(\sigma-1)> = -<(1+(\sigma-1)\tau(\sigma-1))(\sigma-1)(1+(\sigma-1)\tau(\sigma-1))^{-1}, ((\sigma-1)\tau + \tau(\sigma-1))(1+(\sigma-1)\tau(\sigma-1))^{-1}> = -<(1 + \sigma\tau + \sigma^3\tau)(1+\sigma)(1+\sigma\tau+\sigma^3\tau), (\sigma\tau+\sigma^3\tau)(1+\sigma\tau+\sigma^3\tau)> = -<1+\sigma+\sigma^2+\sigma^3, \sigma\tau+\sigma^3\tau> = 0 \in <<i,j> | i, j \in I, <i,j> \in <R,R>>$.

As all mentioned above, we know that: If $<(\sigma-1)^2, \tau-1>$ is non-trival, then it must be the sum of two identical elements of order 4 in $D_1(R)$. To construct the structure of $K_2(R)$, we will introduce the tool of homology of a group:

(Group Module) Let $G$ be an arbitrary group. A $G$-module is an abelian group $M$, together with a map: $G \times M \to M: (g,m) \to gm$, such that, for all $g, g' \in G; m, m' \in M$, we have:
(a) $g(m + m') = gm + gm'$,
(b) $(gg')(m) = g(g')m, 1m = m$.

(Homology of a group) Let $G$ be an arbitrary group. For a $G$-module $M$, define: $M_G = M/\{gm - m | \forall g \in G, \forall m \in M\}$. Now, the functor $Mod_G \to Ab: M \to M^G$ is right exact. Because of that the category of $G$-modules $Mod_G$ has enough projectives, so we have a projective resolution of $M$: $\ldots \to P_2 \xrightarrow{d_2} P_1 \xrightarrow{d_1} P_0 \xrightarrow{d_0} M \to 0$, the complex $\ldots \to (P_2)_G \xrightarrow{d_2} (P_1)_G \xrightarrow{d_1} (P_0)_G \to 0$ need no longer be exact, and we set $H_r(G, M) \stackrel{def}{=} Ker(d_r)/Im(d_{r+1})$ be the $r$-th homology group of $G$-module $M$. Specially, we have: $H_0(G, M) = M_G$.

Lemma 3 For any arbitrary $G$-module $M$, we all have: $H_k(G, M) \cong Tor_k^G(Z, M)$.
Proof of Lemma 3: Let $(Q.) \to M$ is a projective resolution of $M$, then have: $Tor_k^G(Z, M) \stackrel{def}{=}$

$H_k(Z\otimes_{Z[G]}Q.)$ (1), with using the short exact sequence: $0 \to I_G \to Z[G] \to Z \to 0$, the last term $Z$ can be seen as a right $Z[G]$-module (Consider $G$ (right) acts on $Z$ trivially), due to we have: $H_k(G,M) \stackrel{def}{=} H_k((Q.)_G)$ (2), and for any arbitrary $G$-module $N$, we all have: $N_G \stackrel{def}{=} N/\{gn-n\} \cong Z\otimes_{Z[G]}N$, Let $N = Q_0, Q_1, Q_2...$, by comparing with the rigt sides of (1) and (2), we have that: $H_k(G,M) \cong Tor_k^G(Z,M)$, the proof of Lemma3 is completed.

**Lemma 4** For any arbitrary ring $R$, if $n \geq 3$, then have: $(E_n(R))^{ab} \cong 0$.
Proof of Lemma 4: See Lemma 1.3.2 in reference [7].

Back to the original problem, let $M = I_G$, according to Lemma 3, immediately we obtain: $H_1(G, I_G) \cong Tor_1^G(Z, I_G)$, Let the derived functor $Tor_n^G$ acts on the short exact sequence: $0 \to I_G \to Z[G] \to Z \to 0$, we obtain a long exact sequence: $\ldots \to Tor_1^G(Z[G], I_G) \to Tor_1^G(Z, I_G) \stackrel{\partial}{\to} Tor_0^G(I_G, I_G) \to Tor_0^G(Z[G], I_G) \to Tor_0^G(Z, I_G) \to 0$. Because $Z[G]$ is a projective $G$-module, then according to Page 92, Remark A.12 in reference [11], we know that: $Tor_1^G(Z[G], I_G) = 0$. And according to definition, we also can get: $Tor_0^G(I_G, I_G) \cong I_G\otimes_{Z[G]}I_G, Tor_0^G(Z[G], I_G) \cong Z[G]\otimes_{Z[G]}I_G \cong I_G, Tor_0^G(Z, I_G) \cong Z\otimes_{Z[G]}I_G \cong H_0(G, I_G) \cong H_1(G, Z) \cong G^{ab} \stackrel{Lemma\ 4}{\cong} 0$. So as all mentioned above, we get a short exact sequence: $0 \to Tor_1^G(Z, I_G) \to I_G\otimes_{Z[G]}I_G \to I_G \to 0$. And notice that: any an arbitrary torsion element in $I_G\otimes_{Z[G]}I_G$: $\sum x\otimes_{Z[G]}y; x, y \in I_G$, it's image of the map $I_G\otimes_{Z[G]}I_G \to Z[G]\otimes_{Z[G]}I_G \cong I_G \subseteq Z[G]$ in each terms are: $\sum x\otimes_{Z[G]}y \to \sum xy\otimes_{Z[G]}1 \to \sum xy \in Z[G]$, with respectively. Due to $x\otimes_{Z[G]}y$ is a torsion element in $I_G\otimes_{Z[G]}I_G$, so $xy$ is a torsion element in $Z[G]$, and $Z[G]$ is a free $G$-module, so we have: $\sum xy = 0 \Rightarrow : \sum xy\otimes_{Z[G]}1 = 0$, then have: $\sum x\otimes_{Z[G]}y \in Ker(I_G\otimes_{Z[G]}I_G \to Z[G]\otimes_{Z[G]}I_G) = K_2(R)$, i.e.: All elements in $K_2(R)$ are just all torsion elements in $I_G\otimes_{Z[G]}I_G$. Choose an arbitrary non-trivial torsion element in $I_G\otimes_{Z[G]}I_G$: $\sum_{i=1}^{N}(g_i - 1)\otimes_{Z[G]}(S_i - T_i)$, $S_i, T_i \in Z[G]$, and $S_i, T_i$ in each terms are the simplest coefficients combination of elements in $G$, also satisfying the sums of coefficients are equal, the number of terms are permuted as small to large, notice that: $g = e_1e_2...e_w$ where $e_1 \sim e_w$ are generators of the elementary matrices group $G$, i.e.: All in the form: $e_{ij}(x), x \in R$, the order of them are all 2. Then, have: $(g-1)\otimes_{Z[G]}(S-T) = [(e_1 - 1)e_2e_3...e_w + (e_2e_3...e_w - 1)]\otimes_{Z[G]}(S-T) = (e_1 - 1)\otimes_{Z[G]}e_2...e_w(S-T) + (e_2...e_w - 1)\otimes_{Z[G]}(S-T) = \ldots = (e_1 - 1)\otimes_{Z[G]}e_2e_3...e_w(S-T) + (e_2 - 1)\otimes_{Z[G]}e_3...e_w(S-T) + .. + (e_w - 1)\otimes_{Z[G]}(S-T)$, so non-trivial elements in $I_G\otimes_{Z[G]}I_G$ are all in the form: $\sum_{i=1}^{N}(g_i - 1)\otimes_{Z[G]}(S_i - T_i)$, where each $g_i$ is an elementary matrix, and: $ord(g_i) = 2$. According to Def 1.7 in reference [9], we know: $D_1(R)$ is a normal group of $St(2,R)$, as facts we mentioned above, we consider a finite group acts on its normal subgroup by conjugation, then we can define a canonical map sends the homology group to its quotient group as follow:

$\Lambda: H_2(D_1(R), Z) \to H_0(St(2,R)/D_1(R), H_2(D_1(R), Z)) = H_2(D_1(R), Z)/\{sws^{-1} - w, s \in St(2,R), w \in D_1(R)\}: d_1\otimes_{Z[D_1(R)]}d_2 \to \overline{d_1\otimes_{Z[D_1(R)]}d_2}; d_1, d_2 \in I_{D_1(R)}$, where the conjugate action on homology classes is: $d_1\otimes_{Z[D_1(R)]}d_2 \to sd_1s^{-1}\otimes_{Z[D_1(R)]}sd_2s^{-1}$, in fact, by Hochschild-Serre sequence (see reference [5] and [12]), we have a long exact sequence:

$H_2(D_1(R), Z) \stackrel{\phi \circ \Lambda}{\longrightarrow} H_2([R^*, R^*], Z) \to K_2(R) \to H_1(D_1(R), Z) \to H_1([R^*, R^*], Z) \to 0$, $\phi \stackrel{def}{=}$

$$Res_{I_{[R^*,R^*]} \otimes Z_{[R^*,R^*]} I_{[R^*,R^*]}}^{I_{(E_{11}([R^*,R^*]) \oplus E_{22}([R^*,R^*])) \otimes Z_{[E_{11}([R^*,R^*]) \oplus E_{22}([R^*,R^*])]} I_{(E_{11}([R^*,R^*]) \oplus E_{22}([R^*,R^*]))}} \diamond (St(2,R) \to E(2,R))$$ , by

caculation we stated above, we have already know that: $[R^*,R^*] \cong Z_2^{\oplus 3}$, substitute this result into the long exact sequence mentioned above, we immediately get a new long exact sequence:

$H_2(D_1(R), Z) \xrightarrow{\phi \circ \Lambda} H_2(Z_2 \oplus Z_2 \oplus Z_2, Z) \to K_2(R) \to H_1(D_1(R), Z) \to Z_2^{\oplus 3} \to 0$. Next, we will use the Kunneth formula of the finite abelian group version to discuss all cases by classfication. Let's calculate some simple examples first:

In fact: $H_p(Z_2, Z) = \begin{cases} Z_2 & p = 1,3,5,\ldots \\ 0 & p = 2,4,6,\ldots \\ Z & p = 0 \end{cases}$ 及 $H_k = \begin{cases} Z & k = 0 \\ Z_4 & k = 1,3,5,\ldots \\ 0 & k = 2,4,6,\ldots \end{cases}$, so by Kunneth

formula, we have:

$H_1(Z_2 \oplus Z_2, Z) = Z_2 \oplus Z_2, H_2(Z_2 \oplus Z_2, Z) = \oplus_{p=0}^{2}(H_p(Z_2,Z) \otimes_Z H_{2-p}(Z_2,Z)) = (H_0(Z_2,Z) \otimes_Z H_2(Z_2,Z)) \oplus (H_1(Z_2,Z) \otimes_Z H_1(Z_2,Z)) \oplus (H_2(Z_2,Z) \otimes_Z H_0(Z_2,Z)) = Z_2 \otimes_Z Z_2 = Z_2, H_1(Z_2 \oplus Z_2 \oplus Z_2, Z) = Z_2 \oplus Z_2 \oplus Z_2, H_2(Z_2 \oplus Z_2 \oplus Z_2, Z) = (Z \otimes_Z 0) \oplus ((Z_2 \oplus Z_2) \otimes_Z Z_2) \oplus (Z_2 \otimes_Z Z) = Z_2 \oplus Z_2 \oplus Z_2, H_1(Z_4 \oplus Z_4, Z) = Z_4 \oplus Z_4, H_2(Z_4 \oplus Z_4, Z) = Z_4 \otimes_Z Z_4 = Z_4$. Now, we will discuss it in 3 cases:

Case 1. $K_2(R) \cong 0, D_1(R) = Z_2 \oplus Z_2 \oplus Z_2$. Then: $H_2(D_1(R), Z) = Z_2 \oplus Z_2 \oplus Z_2$.

Case2. $K_2(R) \cong Z_2, D_1(R) = Z_2^{\oplus 4}, Z_4 \oplus Z_2 \oplus Z_2$. When: $D_1(R) = Z_2^{\oplus 4}$, have: $H_2(D_1(R), Z) = (Z \otimes_Z Z_2) \oplus (Z_2^{\oplus 4}) \oplus (Z_2 \otimes_Z Z) = Z_2^{\oplus 6}$. When: $D_1(R) = Z_4 \oplus Z_2 \oplus Z_2$ ,have: $H_2(D_1(R), Z) = (Z \otimes_Z Z_2) \oplus (Z_4 \otimes_Z (Z_2 \oplus Z_2)) = Z_2 \oplus Z_2 \oplus Z_2$.

Case 3. $K_2(R) \cong Z_2 \oplus Z_2, D_1(R) = Z_2^{\oplus 5}, Z_4 \oplus Z_2^{\oplus 3}, Z_4 \oplus Z_4 \oplus Z_2$. When: $D_1(R) = Z_2^{\oplus 5}$, we have: $H_2(D_1(R), Z) = (Z \otimes_Z Z_2) \oplus ((Z_2 \oplus Z_2 \oplus Z_2) \otimes_Z (Z_2 \oplus Z_2)) \oplus ((Z_2 \oplus Z_2 \oplus Z_2) \otimes_Z Z) = Z_2^{10}$. When: $D_1(R) = Z_4 \oplus Z_2^{\oplus 3}$, we have: $H_2(D_1(R), Z) = (Z \otimes_Z (Z_2 \oplus Z_2 \oplus Z_2)) \oplus (Z_4 \otimes_Z (Z_2 \oplus Z_2 \oplus Z_2)) = Z_2^{\oplus 6}$. When: $D_1(R) = Z_4 \oplus Z_4 \oplus Z_2$, we have: $H_2(D_1(R), Z) = ((Z_4 \oplus Z_4) \otimes_Z Z_2) \oplus (Z_4 \otimes_Z Z) = Z_4 \oplus Z_2 \oplus Z_2$.

Now, back to the original problem, combining with the definition of $\phi$ : $\phi \stackrel{\text{def}}{=}$

$$Res_{I_{[R^*,R^*]} \otimes Z_{[R^*,R^*]} I_{[R^*,R^*]}}^{I_{(E_{11}([R^*,R^*]) \oplus E_{22}([R^*,R^*])) \otimes Z_{[E_{11}([R^*,R^*]) \oplus E_{22}([R^*,R^*])]} I_{(E_{11}([R^*,R^*]) \oplus E_{22}([R^*,R^*]))}} \diamond (St(2,R) \to E(2,R))$$ , we

know that: If $K_2(R) \cong 0$, then we have: $D_1(R) = K_2(2,R), \phi \circ \Lambda(D_1(R)) = 0$, it is contradicted to the surjection of $\phi \circ \Lambda$, so case 1 can be cancelled. Also notice that: $\forall t x_{12}(a) t^{-1} x_{12}(ua), t x_{21}(a) t^{-1} x_{21}(au^{-1}) \in D_1(R) \cap K_2(2,R); t \mapsto \begin{pmatrix} u & 0 \\ 0 & 1 \end{pmatrix}, a \in R, u \in R^*$ ,

because $[R^*, R^*] \subseteq (\sigma-1)^2 R + 1 = Cent(R^*)$, so if $D_1(R) \cap K_2(2,R) = \{1\}$, then $\phi \circ \Lambda$ is an isomorphism, so there is only $K_2(R) \cong Z_2, D_1(R) = Z_4 \oplus Z_2 \oplus Z_2$ satisfys the condition. If: $D_1(R) \cap K_2(2,R) \supsetneq \{1\}$, i.e.: There exists at least 2 non-trival elements which are not equal: $t x_{12}(a) t^{-1} x_{12}(ua), t x_{21}(a) t^{-1} x_{21}(au^{-1}) \in D_1(R) \cap K_2(2,R)$ , so: $|H_2(D_1(R), Z)| \geq 2^2 \times 2^3 = 2^5$, so only : $K_2(R) \cong Z_2, D_1(R) = Z_2^{\oplus 4}$ or $K_2(R) \cong Z_2 \oplus Z_2, D_1(R) = Z_2^{\oplus 5}$ or $K_2(R) \cong Z_2 \oplus Z_2, D_1(R) = Z_4 \oplus Z_2^{\oplus 3}$ satisfy the condition, then notice we have already know that: If $<(\sigma-1)^2, \tau-1>$ is non-trivial, then it must be the sum of two identical elements of order 4 in $D_1(R)$, so $K_2(R) \cong Z_2 \oplus Z_2, D_1(R) = Z_2^{\oplus 5}$ can be cancelled. Last, notice if: $<(\sigma-1)^2, \tau-1> = 2<\sigma+1, \sigma\tau+\tau> \neq 0$ , and due to: $<\sigma+1, \sigma\tau+\tau> - <\sigma+1, \tau> (By\ R4.) = <\sigma+1, \sigma\tau+(\sigma+1)^2\sigma>$, so: $<\sigma+1, \sigma\tau+\tau> (modI) = <\sigma+1, \tau> (modI) + <\sigma+1, \sigma\tau> (modI)$, it concludes that: $<\sigma+1, \tau>$ is not equal to $<\sigma+1, \sigma\tau+\tau>$, and the orders of them are both

4, so $K_2(R) \cong Z_2 \oplus Z_2, D_1(R) = Z_4 \oplus Z_2^{\oplus 3}$ is also can be cancelled. So as all mentioned above, it can only be: $K_2(R) \cong Z_2$, the proof is finished. ∎

# References

[1] John Milnor, Introduction to Algebraic K-theory, Princeton University Press and University of ToKyo Press, 1971
[2] Bruce A. Magurn, Explicit $K_2$ of some finite group rings, Journal of Pure and Applied Algebra 209 (2007) 801-811
[3] A.W.M.Dress, A.O.KuKu, The Cartan map for equivariant higher algebraic K-groups, Comm. Algebra 9 (1981) 727-746
[4] Jean-Pierre Serre, Local Fields (GTM67), Springer, 1979
[5] 彭喻振，关于稳定秩 1 环的$K_2$群，中国科技大学博士学位论文（2009）
[6] Manfred Kolster, $K_2$ of Non-commutative Local Rings, Journal of Algebra 95, 173-200 (1985)
[7] Charles A.Weibel, The K-book An introduction to Algebraic K-theory (GSM145)
[8] R.K.Dennis, M.R.Stein, $K_2$ of discrete valuation rings
[9] Robert B. Russell, On the relative K2 of non commutative local rings, 1982
[10] M.Kolster, On injective stability for K2. Algebraic K-theory, PartI (Oberwolfach 1980), pp.128-168. Lecture Notes in Mathematics, Vol.966, Springer, Berlin-New York, 1982.
[11] J.S.Milne, Class Fields Theory, http://www.jmilne.org/math.
[12] M.Kolster, General symbols and presentations of elementary linear groups. Journal for die Reine und Angewandte Mathematik, 353 (1984), 132-164.